\documentclass[12pt,reqno]{amsart}
\usepackage{mathrsfs}
\usepackage{latexsym,amsmath,amsfonts,amssymb,amsthm}
 \theoremstyle{plain}
\newtheorem{Thm}{Theorem}
\newtheorem{Lem}{Lemma}

\newtheorem{Conj}{Conjecture}
\newtheorem*{Conj*}{Conjecture}
\theoremstyle{definition}
\newtheorem*{Ack}{Acknowledgment}
\theoremstyle{remark}
\newtheorem*{Rem}{Remark}

\def\A{\mathscr{A}}
\def\S{\mathcal{S}}
\def\I{\mathcal{I}}
\def\J{\mathcal{J}}
\def\F{\mathscr{F}}
\def\L{\mathscr{L}}

\def\R{\mathbb R}

\def\pmod #1{\ ({\rm mod}\ #1)}

\def\mes{{\rm mes}}
\def\floor #1{\lfloor{#1}\rfloor}

\begin{document}
\title{Primes in the form $\floor{\alpha p+\beta}$}
\author{Hongze Li}
\email{lihz@sjtu.edu.cn}
\author{Hao Pan}
\email{haopan79@yahoo.com.cn}
\address{
Department of Mathematics, Shanghai Jiaotong University, Shanghai
200240, People's Republic of China}
\thanks{This work was supported by
the National Natural Science Foundation of China (Grant No.
10771135).} 
\date{}
\maketitle

\begin{abstract}
Let $\beta$ be a real number. Then for almost all irrational
$\alpha>0$ (in the sense of Lebesgue measure)
$$
\limsup_{x\to\infty}\pi_{\alpha,\beta}^*(x)(\log x)^2/x\geq 1,
$$
where
$$
\pi_{\alpha,\beta}^*(x)=\{p\leq x:\, \text{both }p\text{ and
}\floor{\alpha p+\beta}\text{ are primes}\}.
$$
\end{abstract}

Recently Jia \cite{Jia06} solved a conjecture of Long and showed
that for any irrational number $\alpha>0$, there exist infinitely
many primes not in the form $2n+2\floor{\alpha n}+1$, where
$\floor{x}$ denotes the largest integer not exceeding $x$.
Subsequently, in \cite{BanksShparlinski} Banks and Shparlinski
investigated the distribution of primes in the Beatty sequence
$\{\floor{\alpha n+\beta}:\, n\geq 1\}$. Motivated by the binary
Goldbach conjecture and the twin primes conjecture, we have the
following conjecture:
\begin{Conj}
\label{conjecture} Let $\alpha>0$ be an irrational number and
$\beta$ be a real number. Then there exist infinitely many primes
$p$ such that $\floor{\alpha p+\beta}$ is also prime.
\end{Conj}

On the other hand, Deshouillers \cite{Deshouillers76} proved that
for almost all (in the sense of Lebesgue measure) $\gamma>1$ there
exist infinitely many primes $p$ in the form $[n^\gamma]$.
Furthermore, Balog \cite{Balog89} showed that for almost all
$\gamma>1$
$$
\limsup_{x\to\infty}\frac{|\{p\leq x:\, \text{both }p\text{ and
}\floor{p^\gamma}\text{ are primes}\}|}{x/(\log x)^2}\geq\gamma.
$$
In this note we shall show that Conjecture \ref{conjecture} holds
for almost all $\alpha$. Define
$$
\pi_{\alpha,\beta}^*(x)=\{p\leq x:\, \text{both }p\text{ and
}\floor{\alpha p+\beta}\text{ are primes}\}.
$$
\begin{Thm}
\label{theorem} Let $\beta$ be a real number. Then
\begin{equation}
\limsup_{x\to\infty}\pi_{\alpha,\beta}^*(x)(\log x)^2/x\geq 1
\end{equation}
for almost all irrational $\alpha>0$.
\end{Thm}

For a set $X\subseteq\R$, let $\mes(X)$ denote its Lebesgue measure.
Without the additional mentions, the constants implied by $\ll$,
$\gg$ and $O(\cdot)$ will be always absolute.
\begin{Lem}
\label{l1} Let $I\subseteq[0,1)$ be an interval. Suppose that $b,
l>0$. Then
$$
\mes(\{\alpha\in(0,b):\,\{\alpha/l\}\in I\}\})=O(b+l)\mes(I),
$$
where $\{x\}=x-\floor{x}$.
\end{Lem}
\begin{proof}
Without loss of generality, we may assume that $I=(c_1,c_2)$ with
$0\leq c_1<c_2\leq 1$. Let $J=\{\alpha\in(0,b):\,\{\alpha/l\}\in
I\}$. Clearly
\begin{align*}
J\subseteq\bigcup_{0\leq j\leq b/l}((j+c_1)l,\min\{b,(j+c_2)l\}).
\end{align*}
If $l\leq b$, then
$$
\mes(J)\leq(\floor{b/l}+1)(c_2-c_1)l=O(b)\mes(I).
$$
And if $l>b$, then
$$\mes(J)=\mes((c_1l,\min\{b,c_2l\}))\leq
(c_2-c_1)l.
$$
\end{proof}
\begin{Lem}
\label{l2} Suppose that $b_2>b_1>0$ and $\beta$ are arbitrarily
real numbers. Let $\epsilon>0$ be a small number and $x$ be a
sufficiently large (depending on $b_1$, $b_2$, $\beta$ and
$\epsilon$) integer. Then there exists an exceptional set
$J_E\subseteq(b_1,b_2)$ with $\mes(J_E)=O(x^{-\epsilon})$ such
that for any square-free $d\leq x^{1/3-2\epsilon}$ and irrational
$\alpha\in(b_1,b_2)\setminus J_E$,
\begin{equation} |\{1\leq n\leq
x:\,n\floor{\alpha n+\beta}\equiv
0\pmod{d}\}|=\frac{x}{d}\prod_{p\mid
d}\bigg(2-\frac{1}{p}\bigg)+O(x^{1-\epsilon}/d). \end{equation}
\end{Lem}
\begin{proof}
For an irrational $\alpha\in(b_1,b_2)$, let
$$
\A(x;\alpha)=\{n\floor{\alpha n+\beta}:\,1\leq n\leq x\}
$$
and
$$
\A_d(x;\alpha)=\{a\in\A:\,a\equiv 0\pmod{d}\}.
$$
For a square-free $d$, we have
\begin{align*}
|\A_d(x;\alpha)|=&\sum_{s\mid d}|\{1\leq n\leq x/s:\, \floor{\alpha
sn+\beta}\equiv 0\pmod{d/s}, (n,d/s)=1\}|\\
=&\sum_{s\mid d}\sum_{t\mid d/s}\mu(t)|\{1\leq n\leq x/s:\,
\floor{\alpha sn+\beta}\equiv
0\pmod{d/s},t\mid n\}|\\
=&\sum_{s\mid d, t\mid s}\mu(t)|\{1\leq n\leq x/s:\, \floor{\alpha
sn+\beta}\equiv 0\pmod{dt/s}\}|.
\end{align*}
Clearly
\begin{align*}
\floor{\alpha sn+\beta}\equiv
0\pmod{td/s}\Longleftrightarrow\{\alpha ns^2/td+\beta
s/td\}\in[0,s/td).
\end{align*}

Let $\alpha'=\alpha s^2/td$, $\beta'=\beta s/td$, $d'=td/s$ and
$y=x/s$. Clearly $y\geq x^{2/3+2\epsilon}$ and $d'\leq d$. Let
$$
I_{a,q}=\{\theta\in[0,1):\,|\theta q-a|\leq x^{2\epsilon}/y\}.
$$
Suppose that $d'x^{2\epsilon}\leq q\leq y/x^{2\epsilon}$ and $1\leq
a\leq q$ with $(a,q)=1$. If $\{\alpha'\}\in I_{a,q}$,
\begin{align*}
&|\{1\leq n\leq y:\,\{an/q+\beta'\}\in[1/q,1/d'-1/q)\}|\\
\leq&|\{1\leq n\leq y:\,\{\alpha' n+\beta'\}\in[0,1/d')\}|\\
\leq&|\{1\leq n\leq
y:\,\{an/q+\beta'\}\in[0,1/d'+1/q)\cup[1-1/q,1)\}|.
\end{align*}
Hence
$$
|\{1\leq n\leq y:\,\{\alpha'
n+\beta'\}\in[0,1/d')\}|=y/d'+O(q/d')+O(y/q).
$$
Let
$$
\I_{d'}=\bigcup_{\substack{1\leq a\leq q\leq d'x^{2\epsilon}\\
(a,q)=1}}I_{a,q}
$$
Clearly
$$
\mes(\I_{d'})\leq \sum_{\substack{1\leq a\leq q\leq d'x^{2\epsilon}\\
(a,q)=1}}\mes(I_{a,q})\ll\frac{d'x^{4\epsilon}}{y}=
tdx^{4\epsilon-1}.
$$
If $\alpha s^2/td\not\in\I_{td/s}$ for each $s,t$ with $s\mid d,
t\mid s$, then
\begin{align*}
|\A_d(x;\alpha)|=\sum_{s\mid d, t\mid
s}\mu(t)x/td(1+O(x^{-2\epsilon}))=\frac xd\prod_{p\mid
d}(2-1/p)+O(x^{1-\epsilon}/d).
\end{align*}
Let
$$
\J_{d}=\{\alpha\in (0,b):\, \{\alpha s^2/td\}\in\I_{td/s}\text{ for some
}s,t\text{ with }s\mid d, t\mid s\}.
$$
Applying Lemma \ref{l1},
$$
\mes(\J_d)\ll b_2\sum_{\substack{s\mid d, t\mid s\\ b_2\geq
td/s^2}}\mes(\I_{td/s})+\frac{td}{s^2}\sum_{\substack{s\mid d,
t\mid s\\ b_2< td/s^2}}\mes(\I_{td/s})=O(x^{-1/3+\epsilon}).
$$
Finally, Let
$$
J_E=\bigcup_{d\leq x^{1/3-2\epsilon}}\J_d.
$$
Clearly we have $\mes(J_E)=O(x^{-\epsilon})$.
\end{proof}

\begin{Lem}
\label{l3} Suppose that $b_2>b_1>0$, $\epsilon>0$ and $\beta$ are
arbitrarily real numbers. Then there exists an exceptional set
$J_E\subseteq(b_1,b_2)$ with $\mes(J_E)=O(x^{-\epsilon})$ such
that for any irrational $\alpha\in(b_1,b_2)\setminus J_E$,
\begin{equation}
|\{1\leq p\leq x:\,\text{both }p\text{ and }\floor{\alpha
p+\beta}\text{ are primes}\}|\ll\frac{x}{(\log x)^2}
\end{equation}
for sufficiently large (depending on $b_1$, $b_2$,  $\beta$ and
$\epsilon$) $x$.
\end{Lem}
\begin{proof}
Let $z=x^{1/8}$. Define
$$
P(z)=\prod_{\substack{p<z\\ p\text{ prime}}}p
$$
and
$$
\S(A,z)=\{a\in A;\, (a,P(z))=1\}.
$$
Let $\A(\alpha)=\{n\floor{\alpha n+\beta}:\,1\leq n\leq x\}$.
Clearly
$$
\{p\floor{\alpha p+\beta}:\, z+\alpha^{-1}(z+1-\beta)\leq p\leq x,\
\text{both }p\text{ and }\floor{\alpha p+\beta}\text{ are primes}\}
$$
is a subset of $\S(\A(\alpha),z)$. Furthermore, by Lemma \ref{l2},
we know that there exists a set $J_E\subseteq(b_1,b_2)$ with
$\mes(J_E)=O(x^{-\epsilon})$ such that for any square-free $1\leq
d\leq x^{1/3-2\epsilon}$ and irrational
$\alpha\in(b_1,b_2)\setminus J_E$,
$$
|\A_d(\alpha)|=\frac{x}{d}\prod_{p\mid
d}\bigg(2-\frac{1}{p}\bigg)+O(x^{1-\epsilon}/d),
$$
where $\A_d(\alpha)=\{y\in\A(\alpha):\,d\mid y\}$. Let $g(m)$ be
the completely multiplicative function such that $g(p)=2/p-1/p^2$
for each prime $p$. Define $G(z)=\sum_{\substack{m<z\\
m|P(z)}}g(m)$. By Selberg's sieve method,
$$
|\S(\A(\alpha),z)|\leq
\frac{|\A(\alpha)|}{G(z)}+O(\sum_{\substack{d<z^2\\ d\text{
square-free}}}3^{\omega(d)}x^{1-\epsilon}/d),
$$
where $\omega(d)$ denotes the number of distinct prime divisors of
$d$. Since $3^{\omega(d)}\ll d^{\epsilon}$,
$$
\sum_{\substack{d<z^2}}\frac{3^{\omega(d)}}{d}\ll z^{2\epsilon}.
$$
So it suffices to show $G(z)\gg(\log z)^2$. By Theorem 7.14 in
\cite{PanPan92}, we know
$$
G(z)=\sum_{\substack{m<z\\ m\mid
P(z)}}g(m)\gg\prod_{p<z}(1-g(p))^{-1}=\prod_{p<z}(1-2/p+1/p^2)^{-1}\gg(\log
z)^2.
$$
\end{proof}

\begin{proof}[Proof of Theorem \ref{theorem}] Suppose that $b_2>b_1>0$.
Let
$$
\F=\{\alpha\in(b_1,b_2):\,
\limsup_{x\to\infty}\pi_{\alpha,\beta}^*(x)(\log x)^2/x<1\}
$$
and
$$
\F_n=\{\alpha\in(b_1,b_2):\,
\limsup_{x\to\infty}\pi_{\alpha,\beta}^*(x)(\log x)^2/x\leq1-1/n\}.
$$
Clearly $\F=\bigcup_{n>1}\F_n$. So it suffices to show that
$\mes(\F_n)=0$ for every $n>1$. (The measurability of $\F_n$ will be
proven later.)

Assume on the contrary that there exists $n>1$ such that
$\mes(\F_n)>0$. Let $I=(c_1,c_2)$ be an arbitrary sub-interval of
$(b_1,b_2)$. Clearly
\begin{align}
\label{te1} \int_{c_1}^{c_2}\pi_{\alpha,\beta}^*(x)d\alpha=&
\int_{c_1}^{c_2}\bigg(\sum_{\substack{p\leq x\\ p\text{ prime}}}
\sum_{\substack{\alpha p+\beta-1<q\leq\alpha p+\beta\\ q\text{
prime}}}1\bigg)d\alpha\notag\\
=&\sum_{\substack{p\leq x\\ p\text{ prime}}}\sum_{\substack{c_1
p+\beta-1<q\leq c_2 p+\beta\\ q\text{
prime}}}\mes([(q-\beta)/p,(q+1-\beta)/p)\cap[c_1,c_2])\notag\\
\geq&\sum_{\substack{p\leq x\\ p\text{ prime}}}\sum_{\substack{c_1
p+\beta<q\leq c_2 p+\beta-1\\ q\text{ prime}}}\frac{1}{p}\notag\\
=&(c_2-c_1)\sum_{\substack{p\leq x\\
p\text{ prime}}}\frac{1}{\log p}\bigg(1+O\bigg(\frac{1}{\log(c_1
p)}\bigg)\bigg)\notag\\
\geq&(c_2-c_1)\frac{x}{(\log x)^2}\bigg(1+O\bigg(\frac{1}{\log
x}\bigg)\bigg),
\end{align}
provided that $x$ is sufficiently large (depending on $b_1$ and
$b_2$). Suppose that $C>1$ is the implied constant in Lemma
\ref{l3}. Let $\L_I=\F_n\cap I$ and
$$
\L_{I,\delta}(x)=\{\alpha\in I:\,
\pi_{\alpha,\beta}^*(x)\leq(1-\delta)x/(\log x)^2\}.
$$
For any two primes $p$ and $q$, clearly
$$
J_{p,q}:=\{\alpha\in I:\,\floor{\alpha p+\beta}=q\}
$$
is an interval or empty set. Hence
$$
\L_{I,\delta}(x)=I\setminus\bigg(\bigcup_{\substack{k>(1-\delta)x/(\log
x)^2\\ p_1,\ldots,p_k\leq x\text{ are distinct primes}\\
q_1,\ldots,q_k\text{ are primes}}}\bigcap_{j=1}^kJ_{p_j,q_j}\bigg)
$$
is measurable in the sense of Lebesgue measure. Let $\epsilon>0$ be
a very small number. By Lemma \ref{l3},
$$
\int_{c_1}^{c_2}\pi_{\alpha,\beta}^*(x)d\alpha\leq
O(x^{1-\epsilon})+
$$
\begin{equation}
\label{te2} +\mes(\L_{I,\delta}(x))\frac{(1-\delta)x}{(\log
x)^2}+(c_2-c_1-\mes(\L_{I,\delta}(x)))\frac{Cx}{(\log x)^2}
\end{equation}
provided that $x$ is sufficiently large. Combining (\ref{te1}) and
(\ref{te2}), we have
\begin{equation}
\label{te3} \mes(\L_{I,\delta}(x))\leq
\frac{C-1}{C-1+\delta/2}\mes(I).
\end{equation}

We claim that
\begin{equation}
\label{te4} \L_I=\bigcap_{m>n}\bigcup_{y\geq 1}\bigcap_{x\geq
y}\L_{I,1/n-1/m}(x).
\end{equation}
In fact, for any $m>n$, if
$$
\limsup_{x\to\infty}\frac{\pi_{\alpha,\beta}^*(x)}{x/(\log x)^2}<
1-\frac1n+\frac1m,
$$
then there exists $y_0$ such that for any $x\geq y_0$
$$
\pi_{\alpha,\beta}^*(x)\leq\bigg(1-\frac1n+\frac1m\bigg)\frac{x}{(\log
x)^2}.
$$
On the other hand, if $\alpha\in\bigcup_{y}\bigcap_{x\geq
y}\L_{I,1/n-1/m}(x)$, clearly we have
$$
\limsup_{x\to\infty}\frac{\pi_{\alpha,\beta}^*(x)}{x/(\log x)^2}\leq
1-\frac1n+\frac1m.
$$
By (\ref{te3}) and (\ref{te4}), we get
$$
\mes(\L_I)\leq\limsup_{x\to\infty}\L_{I,2/3n}(x)\leq\frac{C-1}{C-1+1/3n}\mes(I).
$$

Since $\mes(\F_n)>0$, there exist open intervals $I_1, I_2,
\ldots\subseteq (b_1,b_2)$ such that
$$
\F_n\subseteq\bigcup_{k=1}^\infty I_k
$$
and
$$
\sum_{k=1}^\infty\mes(I_k)\leq \frac{C-1+1/4n}{C-1}\mes(\F_n).
$$
But by (\ref{te3}),
$$
\mes(\F_n)=\sum_{k=1}^\infty\mes(\L_{I_k})\leq \frac{C-1}{C-1+1/n}
\sum_{k=1}^\infty\mes(I_k)\leq\frac{C-1+1/4n}{C-1+1/3n}\mes(\F_n).
$$
This evidently leads to a contradiction.
\end{proof}

\begin{Rem}
In \cite{Harman88} and \cite{Harman00}, Harman proved that for
almost all real $\alpha>0$ there are infinitely many pairs of
$(p,q)$ satisfying
$$
|\alpha p-q|<\psi(p),\qquad p,q\text{ are primes},
$$
provided that $\psi$ is a non-increasing positive function and
\begin{equation}
\label{diverge}
\sum_{\substack{2\leqslant p\leqslant\infty\\
p\text{ primes}}}\frac{\psi(p)}{\log p}
\end{equation} diverges. (In fact, in \cite{Harman00} Harman established a quantitative version of the above
result, on condition that $\psi(n)\in(0,1/2)$ for each $n$.) As an
immediate consequence, for almost all $\alpha>0$, there exists
infinitely many pair of primes $(p,q)$ such that $[\![\alpha
p]\!]=q$, where $[\![ x]\!]$ is the nearest integer to $x$. For more
related results, the readers may refer to \cite[Chapter
6]{Harman98}.
\end{Rem}

\begin{Ack}
We are grateful to Professor Glyn Harman for his very helpful
discussions and kindly sending us the copies of the references
\cite{Harman88} and \cite{Harman00}.
\end{Ack}

\end{document}